# Dimension reduction for the micromagnetic energy functional on curved thin films

Giovanni Di Fratta

CMAP, École Polytechnique;
route de Saclay, 91128 Palaiseau Cedex - France

Co-affiliation: BCAM - Basque Center for Applied Mathematics;
Alameda de Mazarredo, 14, 48009 Bilbao, Basque Country - Spain.

**Abstract.** Recently a significant interest in ferromagnetic curved thin films has appeared. In particular, thin spherical shells are currently of great interest due to their capability to support skyrmion solutions which can be stabilized by curvature effects only, in contrast to the planar case where the Dzyaloshinskii-Moriya interaction is required. This paper aims to a $\Gamma$-asymptotic analysis of the micromagnetic energy functional, when the shell is generated, like in the case of a sphere, by a bounded, convex and smooth surface.

**A.M.S. Subject Classification:** 82D40, 49S05, 35C20

**Table of contents**

## 1. Introduction and Physical Motivations

### 1.1. The Micromagnetic model

According to the Landau-Lifshitz theory of fine ferromagnetic particles (cf. [8, 9, 11, 12, 30, 38]), the observable states of a rigid ferromagnetic body, occupying a region $\Omega \subseteq \mathbb{R}^3$, are described by its magnetization $M$, a vector field verifying the so-called *fundamental constraint of micromagnetism*: there exists a material dependent positive constant $M_s$ such that $|M| = M_s$ in $\Omega$.

The *spontaneous magnetization* $M_s := M_s(T)$ is highly dependent on the temperature $T$ and vanishes above a critical value $T_c$, characteristic of each crystal type, known as the *Curie point*. Since we will assume the specimen at a fixed temperature well below $T_c$, the value of $M_s$ will be considered constant in $\Omega$. We can, therefore, express the magnetization in the form $M := M_s m$ where $m \colon \Omega \to \mathbb{S}^2$ is a vector field with values in the unit sphere $\mathbb{S}^2$ of $\mathbb{R}^3$.

Although the modulus of $m$ is constant in space, in general, it is not the case for its direction. For single crystal ferromagnets (cf. [1, 5]), the observable states of the magnetization can then be described as the local minimizers of the micromagnetic energy functional which, after a suitable normalization, reads as (cf. [9, p. 22] or [30, p. 138])

$$\mathcal{G}(m,\Omega) := a_{\mathrm{ex}} \underbrace{\int_\Omega |\nabla m|^2}_{=:\mathcal{E}(m)} + \underbrace{\int_\Omega \varphi_{\mathrm{an}}(m)}_{=:\mathcal{A}(m)} - \underbrace{\frac{\mu_0}{2} \int_\Omega \boldsymbol{h}_{\mathrm{d}}[m\chi_\Omega] \cdot m}_{=:\mathcal{W}(m)} - \underbrace{\mu_0 \int_\Omega \boldsymbol{h}_a \cdot m}_{=:\mathcal{Z}(m)}. \qquad (1)$$

with $m \in H^1(\Omega, \mathbb{S}^2)$ and $m\chi_\Omega$ the extension of $m$ by zero outside $\Omega$. The variational analysis of (1) arises as a non-convex and non-local problem.

- The first term, $\mathcal{E}(m)$, is the *exchange energy* and penalizes nonuniformities in the orientation of the magnetization. The positive constant $a_{\mathrm{ex}}$ is the so-called *exchange stiffness constant*, a material specific energy parameter that summarizes the effects of short-range exchange interactions among neighbor spins.

- The *magnetocrystalline anisotropy energy*, $\mathcal{A}(m)$, models the existence of preferred directions of the magnetization (the so-called *easy axes*). The energy density $\varphi_{\mathrm{an}} \colon \mathbb{S}^2 \to \mathbb{R}^+$ is assumed to be a non-negative Lipschitz continuous function that vanishes only on a finite set of directions, the so-called *easy directions*.

- The quantity $\mathcal{W}(m)$ represents the *magnetostatic self-energy* and describes the energy due to the demagnetizing field (stray field) $\boldsymbol{h}_{\mathrm{d}}[m\chi_\Omega]$ generated by $m\chi_\Omega \in L^2(\mathbb{R}^3, \mathbb{R}^3)$. The operator $\boldsymbol{h}_{\mathrm{d}} \colon m \mapsto \boldsymbol{h}_{\mathrm{d}}[m]$ is, for every $m \in L^2(\mathbb{R}^3)$, the unique solution in $L^2(\mathbb{R}^3, \mathbb{R}^3)$ of the Faraday-Maxwell equations of magnetostatics [11, 32] (see section 2.2 for further mathematical details):

$$\begin{cases} \operatorname{div} \boldsymbol{b}[m] = 0, \\ \operatorname{\mathbf{curl}} \boldsymbol{h}_{\mathrm{d}}[m] = \boldsymbol{0}, \\ \boldsymbol{b}[m] = \mu_0(\boldsymbol{h}_{\mathrm{d}}[m] + m), \end{cases} \qquad (2)$$

where $\boldsymbol{b}[m]$ denotes the magnetic flux density and $\mu_0$ is the magnetic permeability of the vacuum.

- Finally, the term $\mathcal{Z}(m)$ is the *interaction energy* and models the tendency of the specimen to have its magnetization aligned with the (externally) applied field $\boldsymbol{h}_a$. The applied field is assumed to be unaffected by variations of $m$.

### 1.2. Physical motivations and state of the art

The four terms in the energy functional (1) take into account effects that originate from different spatial scales, such as short-range exchange forces and long-range magnetostatic interactions. Depending on the relations among the material and geometric parameters of the particle, various asymptotic regimes arise and their investigation can be efficiently carried out by the dimension reduction techniques of calculus of variations (see, e.g., [7, 18, 19, 20, 28, 35, 31, 45]; this list is certainly far from complete).

In this regard, during the last two decades, considerable interest has appeared for magnetic particles having the shape of a curved convex surface (e.g., plane films [19, 28], nanotubes [29, 39, 48], spherical shells [27, 36, 43]). In particular, spherical thin films are currently worth of interest due to their capability to support skyrmion-type solutions which can be stabilized by curvature effects only, in contrast to the planar case where the Dzyaloshinskii-Moriya interaction is required [24, 40]. Nevertheless, and this is the main motivation for our paper, as remarked in [37]: «*It is well established in numerous studies on rigorous micromagnetism that the effects of non-local dipole-dipole interaction can be reduced to an effective easy-surface anisotropy for thin shells when the thickness is much less than the size of the system [..]. Being aware that these results were obtained for plane films, we assume that the same arguments are valid for smoothly curved shells*».

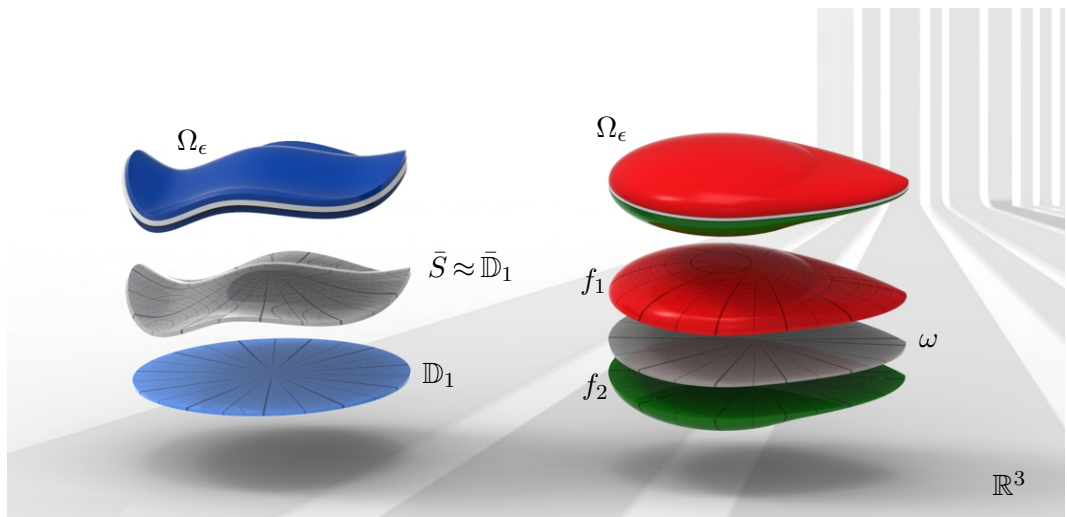

**Figure 1.** (Left) The thin shell $\Omega_\epsilon$ is generated by extruding, in the normal direction $\nu$, a surface $S$ whose closure is globally diffeomorphic to the closed unit disk $\mathbb{D}_1$ of $\mathbb{R}^2$: $\Omega_\epsilon := \{x \in \mathbb{R}^3 : x = \sigma + \epsilon\nu(\sigma), \sigma \in S\}$. (Right) A pillow-like thin shell: $\Omega_\epsilon := \{(x, z) \in \omega \times \mathbb{R} : \epsilon f_2(x) \leqslant z \leqslant \epsilon f_1(x)\}$ where $\omega \subseteq \mathbb{R}^2$ is a planar surface and $f_1$, $f_2$ functions vanishing on the boundary of $\omega$.

Indeed, Gioia and James showed in [28] that for planar thin films the effects of the demagnetizing field operator come down to an effective easy-surface anisotropy. A generalization of this result can be found in [14] where the asymptotic behavior of the energy minimizers is addressed for thin shells generated by surfaces that are diffeomorphic to the closed unit disk of $\mathbb{R}^2$ (see Figure 1). Finally, in [44], a $\Gamma$-convergence analysis is performed on pillow-like shells, i.e., on shells of small thickness $\epsilon > 0$ having the form $\Omega_\epsilon := \{(x, z) \in \omega \times \mathbb{R} : \epsilon f_2(x) \leqslant z \leqslant \epsilon f_1(x)\}$ with $\omega \subseteq \mathbb{R}^2$ and $f_1$, $f_2$ functions vanishing on the boundary of $\omega$ (see Figure 1).

However, all these investigations, being local, do not cover significant physical scenarios like the one of a magnetized thin spherical shell [37]. Nor they can be recovered by a local-to-global gluing argument due to the presence of the *non-local* demagnetizing field operator $\boldsymbol{h}_\mathrm{d}$ (cf. eq. (1)). To have a rigorously justified micromagnetic model of spherical thin films, and more generally of curved convex surfaces, the asymptotic analysis must take into account the global geometry of the surface, and this takes some work. In that respect, the main aim of this paper is the derivation of a reduced model of the micromagnetic energy functional when the region occupied by the ferromagnet is the one of a convex thin film. More specifically, let $S$ be an orientable and smooth convex surface in $\mathbb{R}^3$ (cf. Definition 2), $\nu\colon S\to\mathbb{S}^2$ the normal field associated with the choice of an orientation for $S$. For any sufficiently small $\epsilon>0$ we denote by $\Omega_\epsilon$ the shell having thickness $\epsilon$ and defined by (cf. Figure 2)

$$\Omega_\epsilon := \cup_{\sigma\in S}\ell_\epsilon(\sigma) \quad \text{with} \quad \ell_\epsilon(\sigma):=\{\sigma+t\nu(\sigma)\}_{|t|<\epsilon}. \tag{3}$$

We then consider the energy functional $\mathcal{G}(\cdot,\Omega_\epsilon)$ and use the method of Γ-convergence to characterize the asymptotic behavior of the family $\mathcal{G}(\cdot,\Omega_\epsilon)$ in the limit of vanishing thickness.

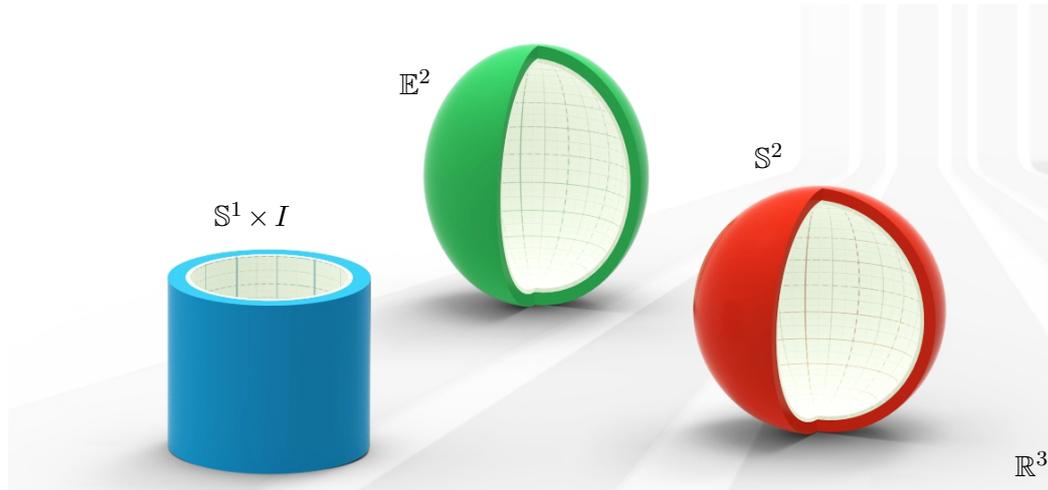

**Figure 2.** Examples of thin shells generated by extruding convex surfaces in their normal direction. (From left to right) The extrusion of a cylinder ($\mathbb{S}^2\times I$), an ellipsoid ($\mathbb{E}^2$) and a sphere ($\mathbb{S}^2$).

The paper is organized as follows: In Section 2 we briefly sketch the geometric setting under which we carry out our investigation. We then state the main result of the paper (cf. Theorem 1) whose proof is given in four steps. The first two steps, developed in Section 3, concern a reformulation of the variational problem and the compactness of minimizing sequences. Finally, Section 4, devoted to the identification of the Γ-limit, completes the proof of the main result.

## 2. Statement of the main result

### 2.1. Notation and setup

We summarize the relevant geometric and functional notions that we use throughout the paper.

#### 2.1.1. Geometric notions

Let $S$ be an orientable and smooth surface in $\mathbb{R}^3$, $\nu\colon S\to\mathbb{S}^2$ the normal field associated with the choice of an orientation for $S$. For every $\sigma\in S$ and every $\delta\in\mathbb{R}^+$ we denote by $\ell_\delta(\sigma):=\{\sigma+t\nu(\sigma)\}_{|t|<\delta}$ the normal segment to $S$ having radius $\delta$ and centered at $\sigma$. We say that $S$ admits a tubular neighborhood if there exists a $\delta\in\mathbb{R}^+$ such that the following properties hold (cf. [23, p. 112]):

**TN$_1$.** For every $\sigma_1,\sigma_2\in S$ one has $\ell_\delta(\sigma_1)\cap\ell_\delta(\sigma_2)=\emptyset$ whenever $\sigma_1\neq\sigma_2$.

**TN$_2$.** The union $\Omega_\delta := \cup_{\sigma \in S} \ell_\delta(\sigma)$ is an open set of $\mathbb{R}^3$ containing $S$. We call $\Omega_\delta$ the tubular neighborhood of $S$ of radius $\delta$.

**TN$_3$.** The map $\psi_\delta \colon (\sigma, s) \in S \times I \mapsto \sigma + s\delta\nu(\sigma) \in \Omega_\delta$, $I := (-1, 1)$, is a diffeomorphism of $S \times I$ onto $\Omega_\delta$. In particular, the nearest point projection $\pi \colon \Omega_\delta \to S$, which maps any $x \in \Omega_\delta$ onto the unique $\sigma \in S$ such that $x \in \ell_\delta(\sigma)$, is a smooth map.

REMARK 1. *Any compact and smooth surface is orientable and admits a tubular neighborhood* [23, Prop. 1, p. 113].

The existence of a tubular neighborhood of $S$ turns out to be *sufficient* to investigate the $\Gamma$-limit of the family of exchange energy functionals $(\mathcal{E}_\epsilon)$ (cf. Proposition 9). On the other hand, due to its non-locality, something more is required for the asymptotic analysis of the family $(\mathcal{W}_\epsilon)$ of magnetostatic self-energies. We, therefore, introduce the following notion.

DEFINITION 2. *Let $S$ be an orientable smooth surface. We say that $S$ is convex if it admits an orientation $\nu$ such that the conditions* **TN$_1$**, **TN$_2$** *and* **TN$_3$**, *still hold when the normal segments $\ell_\delta(\sigma)$ are replaced by the normal half-lines*

$$\ell_\delta^+(\sigma) := \{\sigma + t\nu(\sigma)\}_{t \in (-\delta, +\infty)}. \tag{4}$$

*We then denote by $\Omega_\delta^+$ the unbounded open set $\cup_{\sigma \in S} \ell_\delta^+(\sigma)$ and refer to it as a tubular strip of $S$.*

REMARK 3. Some simple examples of convex surfaces are the sphere $\mathbb{S}^2$ (and more generally the triaxial ellipsoid $\mathbb{E}^2$), the unit cylinder $\mathbb{S}^1 \times I$, the infinite cylinder $\mathbb{S}^1 \times \mathbb{R}$ and the plane $\mathbb{R}^2$ (cf. Figure 2). The name «convex» given to this class of surfaces comes from the compact case where such surfaces are intimately related to the convexity of the domain they bound (cf. [23, Remark 2, p. 393], see also [25, 34]).

For every $\sigma \in S$ the symbols $\tau_1(\sigma), \tau_2(\sigma)$ are used for the orthonormal basis of $T_\sigma S$ made by its principal directions, i.e., the orthonormal basis induced by the eigenvectors of the shape operator of $S$ (cf. [23]). We then write $\kappa_1(\sigma), \kappa_2(\sigma)$ for the principal curvatures at $\sigma \in S$. Note that, when $S$ is convex, the trihedron

$$(\tau_1(\sigma), \tau_2(\sigma), \nu(\sigma)) \quad \text{with} \quad \sigma := \pi(x), x \in \Omega_\delta^+ \tag{5}$$

constitutes a moving frame of $\mathbb{R}^3$ which depends only on $S$. Next, for every $\epsilon \in I_\delta := (0, \delta)$, we introduce the diffeomorphism of $\mathcal{M}$ onto $\Omega_\epsilon$ given by

$$\psi_\epsilon \colon (\sigma, s) \in \mathcal{M} \mapsto \sigma + \epsilon s \nu(\sigma) \in \Omega_\epsilon. \tag{6}$$

Also, we denote by $\mathfrak{g}_\epsilon$ the metric factor which relates the volume form on $\Omega_\epsilon$ to the volume form on $\mathcal{M}$, by $\mathfrak{h}_{1,\epsilon}, \mathfrak{h}_{2,\epsilon}$ the metric coefficients which link the gradient on $\Omega_\epsilon$ to the gradient on $\mathcal{M}$. A direct computation shows that

$$\mathfrak{g}_\epsilon(\sigma, s) := |1 + 2\epsilon s H(\sigma) + (\epsilon s)^2 G(\sigma)|, \quad \mathfrak{h}_{i,\epsilon}(\sigma, s) := \frac{\mathfrak{g}_\epsilon(\sigma, s)}{(1 + \epsilon s \kappa_i(\sigma))^2} \quad (i \in \mathbb{N}_2), \tag{7}$$

where $H(\sigma)$ and $G(\sigma)$ are, respectively, the mean and Gaussian curvature at $\sigma \in S$.

Let $S$ be a smooth and bounded convex surface. We set $I := (-1, 1)$ and denote by $\mathcal{M}$ the product manifold $S \times I$. We then denote by $H^1(\mathcal{M}, \mathbb{R}^3)$ the Sobolev space of vector-valued functions defined on $\mathcal{M}$ (see [2]) endowed with the norm $\|u\|_{H^1(\mathcal{M})}^2 :=$

$\|u\|_{L^2(\mathcal{M})}^2 + \|\nabla_\tau u\|_{L^2(\mathcal{M})}^2 + \|\partial_s u\|_{L^2(\mathcal{M})}^2$ where $\nabla_\tau u$ stands for the tangential gradient of $u$ on $S$. Finally, we write $H^1(\mathcal{M}, \mathbb{S}^2)$ for the subset of $H^1(\mathcal{M}, \mathbb{R}^3)$ made by functions with values in $\mathbb{S}^2$.

2.1.2. THE DEMAGNETIZING FIELD OPERATOR ON $\mathbb{R}^3$

We introduce the Beppo-Levi space

$$W^1(\mathbb{R}^3) = \{u \in \mathcal{S}'(\mathbb{R}^3) : u\omega \in L^2(\mathbb{R}^3), \nabla u \in L^2(\mathbb{R}^3, \mathbb{R}^3)\}, \quad \omega(x) := \frac{1}{\sqrt{1+|x|^2}}, \qquad (8)$$

which is a Hilbert space when endowed with the norm $\|u\|_{W^1(\mathbb{R}^3)}^2 := \|\nabla u\|_{L^2(\mathbb{R}^3, \mathbb{R}^3)}^2$ (cf. [17, Lemma 1, p. 117]. If $\boldsymbol{h}_d[m] \in L^2(\mathbb{R}^3, \mathbb{R}^3)$ is a solution of the Faraday-Maxwell equations (2) then, by Poincaré's lemma [16, Lemma 4, p. 232], there exists a unique magnetostatic potential $u_m \in W^1(\mathbb{R}^3)$ such that $\boldsymbol{h}_d[m] = \nabla u_m$. Hence,

$$-\Delta u_m = \operatorname{div} m \quad \text{in} \quad \mathcal{S}'(\mathbb{R}^3) \qquad (9)$$

and it is straightforward to check, via Lax-Milgram lemma, that for every $m \in L^2(\mathbb{R}^3, \mathbb{R}^3)$ there exists a unique solution of (9) in $W^1(\mathbb{R}^3)$. Therefore, the demagnetizing field can be described as the map which to every magnetization $m \in L^2(\mathbb{R}^3, \mathbb{R}^3)$ associates the distributional gradient of the unique solution of (9) in $W^1(\mathbb{R}^3)$.

REMARK 4. *The weight $\omega$ fix the behavior at infinity of the magnetostatic potential. Note that, in general, $u_m$ does not belong to $L^2(\mathbb{R}^3)$ if $m \in L^2(\mathbb{R}^3, \mathbb{R}^3)$. Indeed, consider any $v \in W^1(\mathbb{R}^3) \setminus H^1(\mathbb{R}^3)$ and set $m = -\nabla v$. We then have $m \in L^2(\mathbb{R}^3, \mathbb{R}^3)$ and $\operatorname{div} m = -\Delta v$. Hence, $u_m := v$ is the unique solution in $W^1(\mathbb{R}^3)$ of (9) and, by construction, $u_m = v \notin L^2(\mathbb{R}^3)$.*

It is easily seen that the map $-\boldsymbol{h}_d : m \in L^2(\mathbb{R}^3, \mathbb{R}^3) \mapsto -\nabla u_m \in L^2(\mathbb{R}^3, \mathbb{R}^3)$ defines a self-adjoint and positive-definite bounded linear operator from $L^2(\mathbb{R}^3, \mathbb{R}^3)$ into itself:

$$-(\boldsymbol{h}_d[m], m)_{L^2(\mathbb{R}^3, \mathbb{R}^3)} = \|\boldsymbol{h}_d[m]\|_{L^2(\mathbb{R}^3, \mathbb{R}^3)}^2 \leqslant \|m\|_{L^2(\mathbb{R}^3, \mathbb{R}^3)}^2 \qquad (10)$$

for every $m \in L^2(\mathbb{R}^3, \mathbb{R}^3)$. Also, notice that for any $\epsilon \in I_\delta$, if $m \in H^1(\Omega_\varepsilon, \mathbb{S}^2)$ then $m\chi_{\Omega_\epsilon} \in L^2(\mathbb{R}^3, \mathbb{R}^3)$ and therefore $\boldsymbol{h}_d[m\chi_{\Omega_\epsilon}] \in L^2(\mathbb{R}^3, \mathbb{R}^3)$. Moreover, since $\boldsymbol{h}_d[m]$ is a gradient field, for any $\epsilon \in I_\delta$ the following two variational equations are satisfied

$$\int_{\mathbb{R}^3} (\boldsymbol{h}_d[m\chi_{\Omega_\epsilon}] + m\chi_{\Omega_\epsilon}) \cdot \nabla \varphi \, \mathrm{d}x = 0, \qquad (11)$$

$$\int_{\mathbb{R}^3} \boldsymbol{h}_d[m\chi_{\Omega_\epsilon}] \cdot \operatorname{\mathbf{curl}} \boldsymbol{\varphi} \, \mathrm{d}x = 0, \qquad (12)$$

for every $\varphi \in W^1(\mathbb{R}^3)$ and any $\boldsymbol{\varphi} \in W^1(\mathbb{R}^3, \mathbb{R}^3)$.

2.1.3. THE ANISOTROPY ENERGY DENSITY AND THE EXTERNAL APPLIED FIELD

The anisotropy energy density $\varphi_{\mathrm{an}} : \mathbb{S}^2 \to \mathbb{R}^+$, which does not depend on $\epsilon \in I_\delta$, is assumed to be a non-negative Lipschitz continuous function that vanishes only on a finite set of directions, the so-called *easy directions*. The hypotheses on $\varphi_{\mathrm{an}}$ are sufficiently general to treat the most common classes of crystal anisotropies arising in applications (e.g., uniaxial, triaxial, cubic). Finally, the external applied field $\boldsymbol{h}_a$ is assumed to be Lipschitz continuous.

2.2. THE MAIN RESULT

To avoid uninformative results (cf. Remark 5) we consider a rescaled version of the energy functional (1). Precisely, let $\Omega_\delta := \cup_{\sigma \in S} \ell_\delta(\sigma)$ be a tubular neighborhood of a smooth convex surface $S$. For any $\epsilon \in I_\delta := (0, \delta)$ we denote by $\mathcal{G}_\epsilon$ the micromagnetic energy functional defined on $H^1(\Omega_\epsilon, \mathbb{S}^2)$ by

$$\mathcal{G}_\epsilon(m) = \frac{1}{\epsilon}(\mathcal{E}_\epsilon(m) + \mathcal{W}_\epsilon(m) + \mathcal{A}_\varepsilon(m) + \mathcal{Z}_\epsilon(m)) \tag{13}$$

$$= \frac{1}{\epsilon}\left(\int_{\Omega_\epsilon} |\nabla m|^2 \, \mathrm{d}x - \frac{1}{2}\int_{\Omega_\epsilon} \boldsymbol{h}_\mathrm{d}[m\chi_{\Omega_\epsilon}] \cdot m \, \mathrm{d}x + \int_{\Omega_\epsilon} \varphi_\mathrm{an}(m) \, \mathrm{d}x - \int_{\Omega_\epsilon} \boldsymbol{h}_a \cdot m \, \mathrm{d}x\right), \tag{14}$$

where $m\chi_{\Omega_\epsilon}$ is the extension of $m$ by zero outside $\Omega_\epsilon$. The existence for any $\epsilon \in I_\delta$ of at least a minimizer for $\mathcal{G}_\epsilon$ in $H^1(\Omega_\epsilon, \mathbb{S}^2)$ is an easy application of the direct method of the calculus of variations (cf. [47]). We are interested in the asymptotic behavior of the family of minimizers of $(\mathcal{G}_\epsilon)_{\epsilon \in I_\delta}$ as $\epsilon \to 0$.

REMARK 5. Let us briefly explain why we rescaled the energy functional (1). It is easily seen that the family $(\epsilon \mathcal{G}_\epsilon)_{\epsilon \in I_\delta}$ is equi-coercive for the weak topology of $H^1(\Omega_\delta)$ and therefore the fundamental theorem of $\Gamma$-convergence applies [10, 15]. On the other hand, as a simple computation shows, the $\Gamma\text{-}\lim_{\epsilon \to 0}(\epsilon \mathcal{G}_\epsilon)$ coincides with the null functional. From this result no information can be retrieved about the asymptotic behavior of the minimizing sequences because every element of $H^1(\Omega_\delta)$ is obviously a minimizer of the null functional.

Next, let us introduce the following functionals defined on $H^1(\mathcal{M}, \mathbb{S}^2)$, which can be thought as the pull-back of $\mathcal{E}_\epsilon, \mathcal{W}_\epsilon, \mathcal{A}_\varepsilon$ and $\mathcal{Z}_\epsilon$ on the product manifold $\mathcal{M} := S \times I$, $I := (-1, 1)$:

- The *exchange energy on* $\mathcal{M}$ reads as $\mathcal{E}^\epsilon_\mathcal{M}(u) := \mathcal{E}^\epsilon_\tau(u) + \mathcal{E}^\epsilon_\nu(u)$, where the tangential and normal component of the exchange energy are respectively given by

$$\mathcal{E}^\epsilon_\tau(u) := \int_\mathcal{M} \sum_{i \in \mathbb{N}_2} |\partial_{\tau_i(\sigma)} u(\sigma, s)|^2 \, \mathfrak{h}_{i,\epsilon}(\sigma, s) \, \mathrm{d}\sigma \, \mathrm{d}s, \tag{15}$$

$$\mathcal{E}^\epsilon_\nu(u) := \frac{1}{\epsilon^2}\int_\mathcal{M} |\partial_s u(\sigma, s)|^2 \, \mathfrak{g}_\epsilon(\sigma, s) \, \mathrm{d}\sigma \, \mathrm{d}s. \tag{16}$$

- The *magnetostatic self-energy on* $\mathcal{M}$ is defined by $\mathcal{W}^\epsilon_\mathcal{M}(u) := \mathcal{W}^\epsilon_\tau(u) + \mathcal{W}^\epsilon_\nu(u)$, where the tangential and normal components of the energy are respectively given by

$$\mathcal{W}^\epsilon_\tau(u) := -\frac{1}{2}\int_\mathcal{M} \sum_{i \in \mathbb{N}_2} (\boldsymbol{h}_\epsilon[u](\sigma, s) \cdot \tau_i(\sigma)) (u(\sigma, s) \cdot \tau_i(\sigma)) \, \mathfrak{g}_\varepsilon(\sigma, s) \, \mathrm{d}\sigma \, \mathrm{d}s, \tag{17}$$

$$\mathcal{W}^\epsilon_\nu(u) := -\frac{1}{2}\int_\mathcal{M} (\boldsymbol{h}_\epsilon[u](\sigma, s) \cdot \nu(\sigma)) (u(\sigma, s) \cdot \nu(\sigma)) \, \mathfrak{g}_\varepsilon(\sigma, s) \, \mathrm{d}\sigma \, \mathrm{d}s. \tag{18}$$

Here, the symbol $\boldsymbol{h}_\epsilon[u] \in L^2(\mathcal{M}, \mathbb{R}^3)$ stands for the *demagnetizing filed on* $\mathcal{M}$:

$$\boldsymbol{h}_\epsilon[u] := \boldsymbol{h}_\mathrm{d}[(u\chi_I) \circ \psi_\epsilon^{-1}] \circ \psi_\epsilon. \tag{19}$$

The family of diffeomorphisms $(\psi_\epsilon)_{\epsilon \in I_\delta}$ is the one given in (6).

- The *anisotropy* and *interaction energies* on $\mathcal{M}$, respectively given by

$$\mathcal{A}_{\mathcal{M}}^{\epsilon}(u) := \int_{\mathcal{M}} \varphi_{\text{an}}(u(\sigma,s))\, \mathfrak{g}_{\varepsilon}(\sigma,s)\, \mathrm{d}\sigma \mathrm{d}s, \qquad (20)$$

$$\mathcal{Z}_{\mathcal{M}}^{\epsilon}(u) := -\int_{\mathcal{M}} \boldsymbol{h}_{a}^{\epsilon}(\sigma,s) \cdot u(\sigma,s)\, \mathfrak{g}_{\varepsilon}(\sigma,s)\, \mathrm{d}\sigma \mathrm{d}s. \qquad (21)$$

For every $\epsilon \in I_\delta$, we have used the symbol $\boldsymbol{h}_a^\epsilon$ for the expression of $\boldsymbol{h}_a$ on $\mathcal{M}$ which, for every $(\sigma,s) \in \mathcal{M}$, is defined by $\boldsymbol{h}_a^\epsilon(\sigma,s) := \boldsymbol{h}_a(\psi_\epsilon(\sigma,s))$. Note that, in the new coordinate system, the applied field depends upon $\epsilon \in I_\delta$.

The main result of the paper is stated in the next result.

THEOREM 1. *For any $\epsilon \in I_\delta$, the minimization problem for $\mathcal{G}_\epsilon$ in $H^1(\Omega_\epsilon, \mathbb{S}^2)$ is equivalent to the minimization in $H^1(\mathcal{M}, \mathbb{S}^2)$ of the functional $\mathcal{F}_\epsilon$ defined by*

$$\mathcal{F}_\epsilon(u) := \mathcal{E}_{\mathcal{M}}^\epsilon(u) + \mathcal{W}_{\mathcal{M}}^\epsilon(u) + \mathcal{A}_{\mathcal{M}}^\epsilon(u) + \mathcal{Z}_{\mathcal{M}}^\epsilon(u), \qquad (22)$$

*in the sense that a configuration $m_\epsilon \in H^1(\Omega_\epsilon, \mathbb{S}^2)$ minimizes $\mathcal{G}_\epsilon$ if and only if $u_\epsilon := m \circ \psi_\epsilon \in H^1(\mathcal{M}, \mathbb{S}^2)$ minimizes $\mathcal{F}_\epsilon$.*

*The family $(\mathcal{F}_\epsilon)_{\epsilon \in I_\delta}$ is equi-coercive in the weak topology of $H^1(\mathcal{M}, \mathbb{S}^2)$ and $\mathcal{F}_0 := \Gamma\text{-}\lim_{\epsilon \to 0} \mathcal{F}_\epsilon$ is given by*

$$\mathcal{F}_0(u) := \mathcal{E}_0(u) + \mathcal{W}_0(u) + \mathcal{A}_0(u) + \mathcal{Z}_0(u) \qquad (23)$$

$$= 2\int_S |\nabla_\tau u|^2 \,\mathrm{d}\sigma + \int_S (u_0 \cdot \nu)^2 \,\mathrm{d}\sigma + 2\int_S \varphi_{\text{an}}(u)\,\mathrm{d}\sigma - 2\int_S \boldsymbol{h}_a \cdot u\,\mathrm{d}\sigma, \qquad (24)$$

*if $\partial_s u = 0$ or by $+\infty$ otherwise. Therefore,*

$$\min_{H^1(\Omega_\epsilon, \mathbb{S}^2)} \mathcal{G}_\epsilon = \min_{H^1(\mathcal{M}, \mathbb{S}^2)} \mathcal{F}_\epsilon = \min_{H^1(\mathcal{M}, \mathbb{S}^2)} \mathcal{F}_0 + o(1) \qquad (25)$$

*and if $(u_\epsilon)_{\epsilon \in I_\delta}$ is a minimizing family for $(\mathcal{F}_\epsilon)_{\epsilon \in I_\delta}$, there exists a subsequence of $(u_\epsilon)_{\epsilon \in I_\delta}$ which weakly converges in $H^1(\mathcal{M}, \mathbb{S}^2)$ to a minimum point of $\mathcal{F}_0$.*

REMARK 6. For the identification of the $\Gamma$-limit of the family of exchange energies, no convexity hypothesis on $S$ is needed. Indeed, following an idea suggested in [14], we only use this assumption to analyze the asymptotic behavior of the magnetostatic self-energy.

REMARK 7. Notice that the families $(\mathcal{A}_{\mathcal{M}}^\epsilon)_{\epsilon \in I_\delta}$ and $(\mathcal{Z}_{\mathcal{M}}^\epsilon)_{\epsilon \in I_\delta}$ constitute a continuous perturbation of $(\mathcal{E}_{\mathcal{M}}^\epsilon + \mathcal{W}_{\mathcal{M}}^\epsilon)_{\epsilon \in I_\delta}$. This means that, with respect to the (topological) product space $I_\delta \times H^1(\mathcal{M}, \mathbb{S}^2)$, the following relations hold

$$\lim_{(\epsilon,u) \to (0,u_0)} \mathcal{A}_{\mathcal{M}}^\epsilon(u) = \int_{\mathcal{M}} \varphi_{\text{an}}(u(\sigma))\,\mathrm{d}\sigma, \qquad (26)$$

$$\lim_{(\epsilon,u) \to (0,u_0)} \mathcal{Z}_{\mathcal{M}}^\epsilon(u) = -\int_{\mathcal{M}} \boldsymbol{h}_a(\sigma) \cdot u(\sigma,s)\,\mathrm{d}\sigma \mathrm{d}s. \qquad (27)$$

Hence, the theorem on the sum of $\Gamma$-limits holds (cf. [15, Prop. 6.20, p. 62]), namely:

$$\Gamma\text{-}\lim_{\epsilon \to 0} \mathcal{F}_\epsilon' = \Gamma\text{-}\lim_{\epsilon \to 0}(\mathcal{E}_{\mathcal{M}}^\epsilon + \mathcal{W}_{\mathcal{M}}^\epsilon) + \Gamma\text{-}\lim_{\epsilon \to 0}\mathcal{A}_{\mathcal{M}}^\epsilon + \Gamma\text{-}\lim_{\epsilon \to 0}\mathcal{Z}_{\mathcal{M}}^\epsilon \qquad (28)$$

$$= \Gamma\text{-}\lim_{\epsilon \to 0}(\mathcal{E}_{\mathcal{M}}^\epsilon + \mathcal{W}_{\mathcal{M}}^\epsilon) + \mathcal{A}_0 + \mathcal{Z}_0. \qquad (29)$$

For this reason, in the identification of the Γ-limit we shall only focus on the family $(\mathcal{E}_{\mathcal{M}}^\epsilon + \mathcal{W}_{\mathcal{M}}^\epsilon)_{\epsilon \in I_\delta}$.

The proof of Theorem 1 is given in four steps. In Subsection 3.1 we prove that for any $\epsilon \in I_\delta$ and any $m \in H^1(\Omega_\epsilon, \mathbb{S}^2)$ the equality $\mathcal{G}_\epsilon(m) = \mathcal{F}_\epsilon(m \circ \psi_\epsilon)$ holds, where $\psi_\epsilon$ stands for the diffeomorphism of $\mathcal{M}$ onto $\Omega_\epsilon$ given by $\psi_\epsilon(\sigma, s) := \sigma + \epsilon s \nu(\sigma)$. In Subsection 3.2 we show that the family $(\mathcal{F}_\epsilon)_{\epsilon \in I_\delta}$ is equi-coercive for the weak topology of $H^1(\mathcal{M}, \mathbb{S}^2)$. Finally, the complete characterization of the Γ-limit $\mathcal{F}_0$ is the object of Section 4.

## 3. Compactness

### 3.1. The equivalence of $\mathcal{G}_\epsilon$ and $\mathcal{F}_\epsilon$

In this section we prove the first part of Theorem 1, namely that once introduced, for any $\epsilon \in I_\delta$, the diffeomorphism of $\mathcal{M}$ onto $\Omega_\epsilon$ given by $\psi_\epsilon: (\sigma, s) \in \mathcal{M} \mapsto \sigma + \epsilon s \nu(\sigma) \in \Omega_\epsilon$, one has $\mathcal{G}_\epsilon(m) = \mathcal{F}_\epsilon(m \circ \psi_\epsilon)$, and therefore $u_\epsilon$ minimizes $\mathcal{G}_\epsilon$ if and only if $u_\epsilon(\sigma, s) := m(\psi_\epsilon(\sigma, \sigma))$ minimizes $\mathcal{F}_\epsilon$.

We only prove the equality $\mathcal{E}_\epsilon(m) = \mathcal{E}_{\mathcal{M}}^\epsilon(m \circ \psi_\epsilon)$, the other ones being easier. For any $m \in H^1(\Omega_\epsilon, \mathbb{S}^2)$, by coarea formula we infer

$$\mathcal{E}_\epsilon(m) := \frac{1}{\epsilon} \int_{\Omega_\epsilon} |\nabla m(x)|^2 \, dx = \int_{\psi_\epsilon(S) \times I} |\nabla m \circ \psi_\epsilon(\sigma, s)|^2 \, d\sigma \, ds \qquad (30)$$

$$= \int_{\mathcal{M}} |\nabla m \circ \psi_\epsilon(\sigma, s)|^2 \, \mathfrak{g}_\varepsilon(\sigma, s) \, d\sigma \, ds. \qquad (31)$$

In writing the last equality we have taken into account that for any $(\epsilon, s) \in I_\delta \times I$ the volume form on $\psi_\varepsilon(S)$ is related to the volume form on $S$ by the metric factor $\mathfrak{g}_\varepsilon(\sigma, s) := |1 + 2\epsilon s H(\sigma) + (\epsilon s)^2 G(\sigma)|$. Next, we project the gradient onto the orthonormal (moving) frame $(\tau_1(\sigma), \tau_2(\sigma), \nu(\sigma))$ induced by $S$ on $\mathbb{R}^3$ (cf. (5)). For any $x \in \Omega_\epsilon$ we have $|\nabla m(x)|^2 = \sum_{i \in \mathbb{N}_2} |\partial_{\tau_i(\sigma)} m(x)|^2 + |\partial_{\nu(\sigma)} m(x)|^2$ with $\sigma = \pi(x)$. Moreover, the following relations hold

$$|\partial_{\tau_i(\sigma)} m(\psi_\epsilon(\sigma, s))|^2 = \frac{1}{(1 + \epsilon s \kappa_i(\sigma))^2} |\partial_{\tau_i(\sigma)} m(\sigma, s)|^2, \qquad (32)$$

$$|\partial_{\nu(\sigma)} m(\psi_\epsilon(\sigma, s))|^2 = \frac{1}{\epsilon^2} |\partial_s m(\sigma, s)|^2, \qquad (33)$$

from which the equality of $\mathcal{E}_\epsilon(m)$ and $\mathcal{E}_{\mathcal{M}}^\epsilon(m \circ \psi_\epsilon)$ results. Note that the previous computation also shows that $m \in H^1(\Omega_\epsilon, \mathbb{S}^2)$ if and only if $m \circ \psi_\epsilon \in H^1(\mathcal{M}, \mathbb{S}^2)$. Finally, as the superposition operator $m \in H^1(\Omega_\epsilon, \mathbb{S}^2) \mapsto (m \circ \psi_\epsilon) \in H^1(\mathcal{M}, \mathbb{S}^2)$ is surjective, we get:

$$\inf_{m \in H^1(\Omega_\epsilon, \mathbb{S}^2)} \mathcal{G}_\epsilon(m) = \inf_{u \in H^1(\mathcal{M}, \mathbb{S}^2)} \mathcal{E}_{\mathcal{M}}^\epsilon(u) + \mathcal{W}_{\mathcal{M}}^\epsilon(u) + \mathcal{A}_{\mathcal{M}}^\epsilon(u) + \mathcal{Z}^\epsilon(u). \qquad (34)$$

This concludes the proof of the first part of Theorem 1.

### 3.2. Equi-coercivity

We now show that the family $(\mathcal{F}_\epsilon)_{\epsilon \in I_\delta}$ is equi-coercive in the weak topology of $H^1(\mathcal{M}, \mathbb{S}^2)$. This means, by definition (see [10]), that there exists a nonempty and weakly compact set $K \subseteq H^1(\mathcal{M}, \mathbb{S}^2)$ such that $\inf_{H^1(\Omega_\epsilon, \mathbb{S}^2)} \mathcal{F}_\epsilon = \inf_K \mathcal{F}_\epsilon$ for every $\epsilon \in I_\delta$.

This is a crucial step in a $\Gamma$-convergence result as it assures the validity of the fundamental theorem of $\Gamma$-convergence concerning the variational convergence of minimum problems ([10, 15]).

Since $\mathcal{A}_\mathcal{M}^\epsilon$ and $\mathcal{Z}_\mathcal{M}^\epsilon$ are uniformly (in $\epsilon \in I_\delta$) bounded terms, it is sufficient to show the equi-coercivity of the family $\mathcal{V}_\epsilon := (\mathcal{E}_\mathcal{M}^\epsilon + \mathcal{W}_\mathcal{M}^\epsilon)_{\epsilon \in I_\delta}$. To this end, we observe that for any constant in space $v \in H^1(\mathcal{M}, \mathbb{S}^2)$ we have

$$\min_{u \in H^1(\mathcal{M}, \mathbb{S}^2)} \mathcal{V}_\epsilon(u) \leqslant \mathcal{E}_\mathcal{M}^\epsilon(v) + \mathcal{W}_\mathcal{M}^\epsilon(v) = \mathcal{W}_\mathcal{M}^\epsilon(v). \tag{35}$$

Taking into account (10) and that $\mathfrak{g}_\epsilon$ is bounded on $\mathcal{M}$, uniformly with respect to $\epsilon \in I_\delta$, we end up with

$$\min_{u \in H^1(\mathcal{M}, \mathbb{S}^2)} \mathcal{V}_\epsilon(u) \leqslant \int_\mathcal{M} \mathfrak{g}_\epsilon(\sigma, s) \, d\sigma ds \leqslant \kappa_\mathcal{M} |\mathcal{M}|, \tag{36}$$

for a suitable positive constant $\kappa_\mathcal{M}$ depending only on $\mathcal{M}$. Therefore, for every $\epsilon \in I_\delta$, the minimizers of $(\mathcal{V}_\epsilon)_{\epsilon \in I_\delta}$ are in $K(\mathcal{M}, \mathbb{S}^2) := \cup_{\epsilon \in I_\delta} \{u \in H^1(\mathcal{M}, \mathbb{S}^2) : \mathcal{V}_\epsilon(u) \leqslant \kappa_\mathcal{M} |\mathcal{M}|\}$. On the other hand, since the principal curvatures $\kappa_1, \kappa_2$ are bounded in $S$, whenever the radius $\delta \in \mathbb{R}^+$ of the tubular neighborhood $\Omega_\delta$ is sufficiently small, there exists a positive constant $c_\mathcal{M}$, independent from $\epsilon \in I_\delta$, such that for any $i \in \mathbb{N}_2$ one has $\inf_{(\sigma, s) \in \mathcal{M}} \mathfrak{h}_{i,\epsilon}(\sigma, s) \geqslant c_\mathcal{M}$ for every $\epsilon \in I_\delta$. Therefore, since $\mathcal{W}_\mathcal{M}^\epsilon$ is always nonnegative because of (10), we get

$$\|u\|_{H^1(\mathcal{M}, \mathbb{S}^2)}^2 = |\mathcal{M}| + \sum_{i \in \mathbb{N}_2} \int_\mathcal{M} |\partial_{\tau_i(\sigma)} u(\sigma, s)|^2 \, d\sigma ds + \int_\mathcal{M} |\partial_s u(\sigma, s)|^2 \, d\sigma ds$$

$$\leqslant |\mathcal{M}| + \frac{1}{c_\mathcal{M}} \mathcal{V}_\epsilon(u), \tag{37}$$

and therefore if $u \in K(\mathcal{M}, \mathbb{S}^2)$ then $\|u\|_{H^1(\mathcal{M}, \mathbb{S}^2)}^2 \leqslant (1 + \kappa_\mathcal{M}/c_\mathcal{M})|\mathcal{M}|$. In other words, the set $K(\mathcal{M}, \mathbb{S}^2)$ is contained in the bounded subset $H_b^1(\mathcal{M}, \mathbb{S}^2)$ of $H^1(\mathcal{M}, \mathbb{R}^3)$ given by the intersection of $H^1(\mathcal{M}, \mathbb{S}^2)$ with the ball of $H^1(\mathcal{M}, \mathbb{R}^3)$ centered at the origin and of radius $1 + \kappa_\mathcal{M}/c_\mathcal{M}$. Thus, for any $\epsilon \in I_\delta$

$$\min_{u \in H^1(\mathcal{M}, \mathbb{S}^2)} \mathcal{V}_\epsilon(u) = \min_{u \in H_b^1(\mathcal{M}, \mathbb{S}^2)} \mathcal{V}_\epsilon(u). \tag{38}$$

To prove that $H_b^1(\mathcal{M}, \mathbb{S}^2)$ is weakly compact it is sufficient to show that the set $H_b^1(\mathcal{M}, \mathbb{S}^2)$ is weakly closed. To this end, we note that if $(u_n)_{n \in \mathbb{N}}$ is a sequence in $H_b^1(\mathcal{M}, \mathbb{S}^2)$ such that $u_n \rightharpoonup u_0$ weakly in $H^1(\mathcal{M}, \mathbb{R}^3)$, due to Rellich-Kondrachov theorem, $u_n \to u_0$ strongly in $L^2(\mathcal{M}, \mathbb{R}^3)$, and therefore, up to the extraction of a subsequence, $1 \equiv |u_n| \to |u_0|$ a.e. in $\mathcal{M}$. Thus $u_0(\sigma, s) \in \mathbb{S}^2$ for a.e. $(\sigma, s) \in \mathcal{M}$ and this concludes the proof.

### 4. The identification of the $\Gamma$-limit

In this section, we compute $\mathcal{F}_0 := \Gamma\text{-}\lim_{\epsilon \to 0} \mathcal{F}_\epsilon$. As pointed out at the end of Subsection 2.2, it is sufficient to focus on the $\Gamma$-convergence of the family

$$\mathcal{V}_\epsilon : u \in H^1(\mathcal{M}, \mathbb{S}^2) \mapsto \mathcal{E}_\mathcal{M}^\epsilon(u) + \mathcal{W}_\mathcal{M}^\epsilon(u). \tag{39}$$

We set $\mathcal{V}_0 := \mathcal{E}_0 + \mathcal{W}_0$ with $\mathcal{E}_0$ and $\mathcal{W}_0$ given by (23). Note that, as a consequence of (10), $\mathcal{V}_\epsilon(u) \geqslant 0$ for any $u \in H^1(\mathcal{M}, \mathbb{S}^2)$.

Let us prove the Γ-liminf inequality for $(\mathcal{V}_\epsilon)_{\epsilon \in I_\delta}$, i.e., that for any family $(u_\epsilon)_{\epsilon \in I_\delta}$ weakly convergent to some $u_0 \in H^1(\mathcal{M}, \mathbb{S}^2)$ we have $\mathcal{V}_0(u_0) \leqslant \liminf_{\epsilon \to 0} \mathcal{V}_\epsilon(u_\epsilon)$. With no loss of generality, we can assume that $\liminf_{\epsilon \to 0} \mathcal{V}_\epsilon(u_\epsilon) < +\infty$. We then have (see (16))

$$+\infty > \liminf_{\epsilon \to 0} \mathcal{F}_\epsilon(u_\epsilon) \geqslant \liminf_{\epsilon \to 0} \mathcal{E}_\nu^\epsilon(u_\epsilon)$$
$$= \liminf_{\epsilon \to 0} \frac{1}{\epsilon^2} \int_{\mathcal{M}} |\partial_s u_\epsilon(\sigma, s)|^2 \mathfrak{g}_\varepsilon(\sigma, s) \, \mathrm{d}\sigma \mathrm{d}s. \quad (40)$$

Moreover, for any $i \in \mathbb{N}_2$, since $\sup_{\sigma \in S} |\kappa_i(\sigma)| < \infty$, there exists a strictly positive real-valued function $\gamma: I_\delta \to \mathbb{R}^+$ such that, at least in a neighborhood of $0 \in \mathbb{R}$, the following estimate holds:

$$\inf_{(\sigma,s) \in \mathcal{M}} \mathfrak{h}_{i,\epsilon}(\sigma, s) = \inf_{(\sigma,s) \in \mathcal{M}} \frac{\mathfrak{g}_\varepsilon(\sigma,s)}{(1+\epsilon s \kappa_i(\sigma))^2} \geqslant \gamma(\epsilon) \quad \text{with} \quad \gamma(\epsilon) = 1 + o(1). \quad (41)$$

Using (40) and (41) we find that $\lim_{\epsilon \to 0} \|\partial_s u_0\|_{L^2(\mathcal{M})} = 0$. Since $\partial_s u_\epsilon \rightharpoonup \partial_s u_0$ in $\mathcal{D}'(\mathcal{M})$ we infer that

$$\partial_s u_\epsilon \to \partial_s u_0(\sigma, s) \text{ strongly in } L^2(\mathcal{M}), \quad \partial_s u_0(\sigma, s) = 0 \text{ a.e. in } \mathcal{M}. \quad (42)$$

Therefore, for the identification of the Γ-limit of $(\mathcal{V}_\epsilon)_{\epsilon \in I_\delta}$ it is sufficient to restrict the analysis to the families of $H^1(\mathcal{M}, \mathbb{S}^2)$ functions which weakly converge to an element $u_0 \in H^1(\mathcal{M}, \mathbb{S}^2)$ having the form

$$u_0(\sigma, s) = \chi_I(s) \tilde{u}_0(\sigma), \quad (43)$$

for some $\tilde{u}_0 \in H^1(S, \mathbb{S}^2)$, i.e., not depending on the $s$ variable. In the following, with a slight abuse of notation, we shall write $u_0(\sigma)$ instead of $\tilde{u}_0(\sigma)$.

In computing $\mathcal{V}_0$, we first show that the Γ-limit of the families $(\mathcal{E}_\mathcal{M}^\epsilon)_{\epsilon \in I_\delta}$ and $(\mathcal{W}_\mathcal{M}^\epsilon)_{\epsilon \in I_\delta}$, is respectively equal to $\mathcal{E}_0$ and $\mathcal{W}_0$ (cf. (23)), then we prove that $\mathcal{V}_0 := \Gamma\text{-}\lim_{\epsilon \to 0} \mathcal{V}_\epsilon = \mathcal{E}_0 + \mathcal{W}_0$.

### 4.1. The Γ-limit of the family $(\mathcal{E}_\mathcal{M}^\epsilon)_{\epsilon \in I_\delta}$

This section provides the identification of the Γ-limit of the family of exchange energies on $\mathcal{M}$. Note that, in what follows, we will *not* make use of the convexity assumption on $S$. We start by addressing the Γ-liminf inequality for $(\mathcal{E}_\mathcal{M}^\epsilon)_{\epsilon \in I_\delta}$. Taking into account the lower semicontinuity of the norm, for any $u_\epsilon \rightharpoonup u_0$ in $H^1(\mathcal{M}, \mathbb{S}^2)$, with $u_0$ of the type (43), we get

$$\|u_0\|_{H^1(\mathcal{M}, \mathbb{S}^2)}^2 = \int_\mathcal{M} |u_0(\sigma, s)|^2 \mathrm{d}\sigma \mathrm{d}s + \sum_{i \in \mathbb{N}_2} \int_\mathcal{M} |\partial_{\tau_i(\sigma)} u_0(\sigma)|^2 \mathrm{d}\sigma \mathrm{d}s \quad (44)$$

$$\leqslant |\mathcal{M}| + \liminf_{\epsilon \to 0} \int_\mathcal{M} \sum_{i \in \mathbb{N}_2} |\partial_{\tau_i(\sigma)} u_\epsilon(\sigma, s)|^2 + |\partial_s u_\epsilon(\sigma, s)|^2 \, \mathrm{d}\sigma \mathrm{d}s \quad (45)$$

$$= |\mathcal{M}| + \liminf_{\epsilon \to 0} \int_\mathcal{M} |\nabla_{\tau(\sigma)} u_\epsilon(\sigma, s)|^2 \mathrm{d}\sigma \mathrm{d}s. \quad (46)$$

In deriving the last equality, we used (42) and denoted by $\nabla_{\tau(\sigma)} u_\epsilon(\sigma, s)$ the tangential gradient of $u_\epsilon$ on $S$ whose norm, with respect to an orthonormal basis $(\tau_1(\sigma), \tau_2(\sigma))$ of $T_\sigma S$, can be expressed as $|\nabla_{\tau(\sigma)} u_\epsilon(\sigma, s)|^2 := \sum_{i \in \mathbb{N}_2} |\partial_{\tau_i(\sigma)} u_\epsilon(\sigma, s)|^2$. Thus

$$2\|\nabla_\tau u_0\|_{H^1(S, \mathbb{S}^2)}^2 = 2 \sum_{i \in \mathbb{N}_2} \int_S |\partial_{\tau_i(\sigma)} u_0(\sigma)|^2 \mathrm{d}\sigma \mathrm{d}s \leqslant \liminf_{\epsilon \to 0} \int_\mathcal{M} |\nabla_\tau u_\epsilon|^2. \quad (47)$$

Next, by making use of the well-known properties of the liminf operator and taking into account relation (41), we compute:

$$\liminf_{\epsilon\to 0}\int_{\mathcal{M}}|\nabla_\tau u_\epsilon|^2 = \left(\liminf_{\epsilon\to 0}\gamma(\epsilon)\right)\left(\liminf_{\epsilon\to 0}\int_{\mathcal{M}}|\nabla_{\tau(\sigma)}u_\epsilon(\sigma,s)|^2\,\mathrm{d}\sigma\,\mathrm{d}s\right) \quad (48)$$

$$\leqslant \liminf_{\epsilon\to 0}\left(\gamma(\epsilon)\int_{\mathcal{M}}|\nabla_{\tau(\sigma)}u_\epsilon(\sigma,s)|^2\,\mathrm{d}\sigma\,\mathrm{d}s\right) \quad (49)$$

$$= \liminf_{\epsilon\to 0}\left(\sum_{i\in\mathbb{N}_2}\gamma(\epsilon)\int_{\mathcal{M}}|\partial_{\tau_i(\sigma)}u_\epsilon(\sigma,s)|^2\,\mathrm{d}\sigma\,\mathrm{d}s\right) \quad (50)$$

$$\leqslant \liminf_{\epsilon\to 0}\left(\sum_{i\in\mathbb{N}_2}\int_{\mathcal{M}}|\partial_{\tau_i(\sigma)}u_\epsilon(\sigma,s)|^2\frac{\mathfrak{g}_\varepsilon(\sigma,s)}{(1+\epsilon s\kappa_i(\sigma))^2}\,\mathrm{d}\sigma\,\mathrm{d}s\right) \quad (51)$$

$$= \liminf_{\epsilon\to 0}\mathcal{E}_\tau^\epsilon(u_\epsilon). \quad (52)$$

Substituting (51) into (47) we get the following result.

LEMMA 8. *Suppose that $u_\epsilon \rightharpoonup u_0$ weakly in $H^1(\mathcal{M},\mathbb{S}^2)$ and $\liminf_{\epsilon\to 0}\mathcal{V}_\epsilon(u_\epsilon) < +\infty$. The following estimate holds*

$$2\|\nabla_\tau u_0\|_{H^1(S,\mathbb{S}^2)}^2 \leqslant \liminf_{\epsilon\to 0}\mathcal{E}_\tau^\epsilon(u_\epsilon) \leqslant \liminf_{\epsilon\to 0}(\mathcal{E}_\tau^\epsilon(u_\epsilon)+\mathcal{E}_\nu^\epsilon(u_\epsilon)). \quad (53)$$

We now address the existence of a recovery sequence. To this end, it is sufficient to note that for every $u_\epsilon \in H^1(\mathcal{M},\mathbb{S}^2)$ having the product form $u_\epsilon(\sigma,s) = \chi_I(s)u_0(\sigma)$ we have (cf. (41))

$$\limsup_{\epsilon\to 0}(\mathcal{E}_\tau^\epsilon(u_\epsilon)+\mathcal{E}_\nu^\epsilon(u_\epsilon)) = \limsup_{\epsilon\to 0}\sum_{i\in\mathbb{N}_2}\int_{\mathcal{M}}|\partial_{\tau_i(\sigma)}u_0(\sigma)|^2\frac{\mathfrak{g}_\varepsilon(\sigma,s)}{(1+\epsilon s\kappa_i(\sigma))^2}\,\mathrm{d}\sigma\,\mathrm{d}s \quad (54)$$

$$= \sum_{i\in\mathbb{N}_2}\int_{\mathcal{M}}|\partial_{\tau_i(\sigma)}u_0(\sigma)|^2\,\mathrm{d}\sigma\,\mathrm{d}s \quad (55)$$

$$= 2\|\nabla_\tau u_0\|_{L^2(S,\mathbb{S}^2)}^2. \quad (56)$$

We so proved the following result.

PROPOSITION 9. *Let $S$ be a smooth compact surface (convex or not) and $\mathcal{M} := S \times I$. The family $(\mathcal{E}_{\mathcal{M}}^\epsilon)_{\epsilon\in I_\delta}$ of exchange energy on $\mathcal{M}$, $\Gamma$-converges, with respect to the weak topology of $H^1(\mathcal{M},\mathbb{S}^2)$, to the functional*

$$\mathcal{E}_0 : u \in H^1(\mathcal{M},\mathbb{S}^2) \mapsto \begin{cases} 2\|\nabla_\tau u\|_{L^2(S,\mathbb{S}^2)}^2 & \text{if } \partial_s u = 0, \\ +\infty & \text{otherwise}. \end{cases} \quad (57)$$

### 4.2. THE $\Gamma$-LIMIT OF THE FAMILY $(\mathcal{W}_{\mathcal{M}}^\epsilon)_{\epsilon\in I_\delta}$

This section is devoted to the identification of the $\Gamma$-limit of the family $(\mathcal{W}_{\mathcal{M}}^\epsilon)_{\epsilon\in I_\delta}$ of magnetostatic self-energies on $\mathcal{M}$.

Note that for every $u \in L^2(\mathcal{M}, \mathbb{R}^3)$ the distribution $(u\chi_I) \circ \psi_\epsilon^{-1}$, with $\psi_\epsilon$ given by (6), is in $L^2(\Omega_\epsilon, \mathbb{R}^3)$, and it is therefore possible to evaluate the demagnetizing field $\boldsymbol{h}_\mathrm{d}$ on its extension by zero outside $\Omega_\epsilon$. To simplify the notation, we still denote by $(u\chi_I) \circ \psi_\epsilon^{-1}$ such an extension. Since $S$ is a convex surface (cf. Definition 2), for every $\epsilon \in I_\delta$ there exists the tubular strip of $S$. Namely

$$\Omega_{\mathcal{M}_+}^\epsilon := \{\sigma + \epsilon s \nu(\sigma)\}_{(\sigma,s) \in \mathcal{M}_+} \tag{58}$$

with $\mathcal{M}_+ := S \times (-1, +\infty)$. We consider the restriction of (11) and (12) to $\Omega_{\mathcal{M}_+}^\epsilon$ and pull-back them via the diffeomorphism $\psi_\epsilon \colon (\sigma, s) \in \mathcal{M}_+ \mapsto \sigma + \epsilon s \nu(\sigma) \in \Omega_{\mathcal{M}_+}$. We so obtain the following relations:

$$\epsilon \int_{\mathcal{M}_+} (\boldsymbol{h}_\epsilon[u](\sigma, s) + u(\sigma, s)\chi_I(s)) \cdot (\nabla\varphi \circ \psi_\epsilon) \, \mathfrak{g}_\epsilon(\sigma, s) \, \mathrm{d}\sigma \mathrm{d}s = 0 \tag{59}$$

$$\epsilon \int_{\mathcal{M}_+} \boldsymbol{h}_\epsilon[u](\sigma, s) \cdot (\mathbf{curl}\boldsymbol{\varphi} \circ \psi_\epsilon) \, \mathfrak{g}_\epsilon(\sigma, s) \, \mathrm{d}\sigma \mathrm{d}s = 0 \tag{60}$$

for any $\varphi \in \mathcal{D}(\Omega_{\mathcal{M}_+}^\epsilon)$, $\boldsymbol{\varphi} \in \mathcal{D}(\Omega_{\mathcal{M}_+}^\epsilon, \mathbb{R}^3)$. Here, $\boldsymbol{h}_\epsilon[u](\sigma, s) := \boldsymbol{h}_\mathrm{d}[(u\chi_I) \circ \psi_\epsilon^{-1}] \circ \psi_\epsilon$.

Next, let $(u_\epsilon)_{\epsilon \in I_\delta}$ be a family of $H^1(\mathcal{M}, \mathbb{S}^2)$ functions weakly converging to some $u_0 \in H^1(\mathcal{M}, \mathbb{S}^2)$ and such that $\liminf_{\epsilon \to 0} \mathcal{F}_\epsilon(u_\epsilon) < +\infty$. By relations (10) and (19), we deduce that for any $m_\epsilon \chi_{\Omega_\epsilon} := (u_\epsilon \chi_I) \circ \psi_\epsilon^{-1}$

$$\frac{1}{2} \int_{\mathcal{M}_+} |\boldsymbol{h}_\epsilon[u_\epsilon](\sigma, s)|^2 \mathfrak{g}_\epsilon(\sigma, s) \, \mathrm{d}\sigma \mathrm{d}s = \frac{1}{2\epsilon} \int_{\Omega_{\mathcal{M}_+}^\epsilon} |\boldsymbol{h}_\mathrm{d}[m_\epsilon \chi_{\Omega_\epsilon}]|^2 \, \mathrm{d}\mu \tag{61}$$

$$\leqslant \frac{1}{2\epsilon} \int_{\mathbb{R}^3} |\boldsymbol{h}_\mathrm{d}[m_\epsilon \chi_{\Omega_\epsilon}]|^2 \, \mathrm{d}\mu \tag{62}$$

$$= -\frac{1}{2\epsilon} \int_{\Omega_\epsilon} \boldsymbol{h}_\mathrm{d}[m_\epsilon \chi_{\Omega_\epsilon}] \cdot m_\epsilon \chi_{\Omega_\epsilon} \, \mathrm{d}\mu \tag{63}$$

$$= \mathcal{W}_\tau^\epsilon(u_\epsilon) + \mathcal{W}_\nu^\epsilon(u_\epsilon), \tag{64}$$

with $\mathcal{W}_\tau^\epsilon$ and $\mathcal{W}_\nu^\epsilon$ respectively given by (17) and (18). Hence, there exist a subsequence extracted from $(\boldsymbol{h}_\epsilon[u_\epsilon\chi_I])_{\epsilon \in I_\delta}$, still denoted by $(\boldsymbol{h}_\epsilon[u_\epsilon\chi_I])_{\epsilon \in I_\delta}$, and an element $\boldsymbol{h}_0 \in L^2(\mathcal{M}_+, \mathbb{R}^3)$, such that $\boldsymbol{h}_\epsilon[u_\epsilon\chi_I] \rightharpoonup \boldsymbol{h}_0$ weakly in $L^2(\mathcal{M}_+, \mathbb{R}^3)$.

Let us consider the energy term $\mathcal{W}_\nu$, i.e., the normal part of the family of magnetostatic self-energy functionals defined by (18). Decomposing (59) into its normal and tangential part, and evaluating it on the weakly convergent sequence $(u_\epsilon)_{\epsilon \in I_\delta}$ we get that

$$\int_{\mathcal{M}_+} [(\boldsymbol{h}_\epsilon[u_\epsilon](\sigma, s) + u_\epsilon(\sigma, s)\chi_I(s)) \cdot \nu(\sigma)]\partial_s\varphi(\sigma, s) \, \mathfrak{g}_\epsilon(\sigma, s) \, \mathrm{d}\sigma \mathrm{d}s$$
$$= -\epsilon \sum_{i \in \mathbb{N}_2} \int_{\mathcal{M}^+} (\boldsymbol{h}_\epsilon[u_\epsilon](\sigma, s) + u_\epsilon(\sigma, s)\chi_I(s)) \cdot \tau_i(\sigma) \, \partial_{\tau_i}\varphi(\sigma, s) \, \mathfrak{h}_{i,\epsilon}(\sigma, s) \, \mathrm{d}\sigma \mathrm{d}s. \tag{65}$$

for any $\varphi \in \mathcal{D}(\mathcal{M}_+)$. Taking into account (41) and passing to the limit for $\epsilon \to 0$ in (65) we get, up to the extraction of a subsequence,

$$\int_{\mathcal{M}_+} [(\boldsymbol{h}_0(\sigma, s) + u_0(\sigma)\chi_I(s)) \cdot \nu(\sigma)]\partial_s\varphi(\sigma, s) \, \mathrm{d}\sigma \mathrm{d}s = 0 \tag{66}$$

for any $\varphi \in \mathcal{D}(\mathcal{M}_+)$. Thus the quantity $(\boldsymbol{h}_0(\sigma, s) + u_0(\sigma)\chi_I(s)) \cdot \nu(\sigma)$ is constant with respect to the $s$-variable, and since it belongs to $L^2(\mathcal{M}_+)$ we infer that

$$\boldsymbol{h}_0(\sigma, s) \cdot \nu(\sigma) = -u_0(\sigma)\chi_I(s) \cdot \nu(\sigma) \tag{67}$$

for a.e. $(\sigma, s) \in \mathcal{M}_+$. In particular, the normal component of the weak limit $\boldsymbol{h}_0 \in L^2(\mathcal{M}_+, \mathbb{R}^3)$ does not depend on the extracted subsequence so that the full subsequence $\boldsymbol{h}_\epsilon[u_\epsilon \chi_I](\sigma, s) \cdot \nu(\sigma)$ weakly converges to $-u_0(\sigma) \chi_I(s) \cdot \nu(\sigma)$ in $L^2(\mathcal{M}_+, \mathbb{R}^3)$. By Rellich–Kondrachov theorem, the weak convergence of $(u_\epsilon)_{\epsilon \in I_\delta}$ to $u_0(\sigma) \in H^1(\mathcal{M}, \mathbb{S}^2)$ implies that $u_\epsilon(\sigma, s) \to u_0(\sigma)$ strongly in $L^2(\mathcal{M}, \mathbb{R}^3)$. By taking the limit for $\epsilon \to 0$ of both members of (65), taking into account (67), we finish with the following relation:

$$\lim_{\epsilon \to 0} \mathcal{W}_\nu^\epsilon(u_\epsilon) = \frac{1}{2} \int_{\mathcal{M}} (u_0(\sigma) \cdot \nu(\sigma))^2 \, d\sigma \, ds. \tag{68}$$

Note that the right-hand side of (68) coincides with $\mathcal{W}_0$ because, as we are going to show, the demagnetizing field $\boldsymbol{h}_0$ has no tangential component.

We now address the tangential energy term $\mathcal{W}_\tau^\epsilon$ defined by (17). We start by decomposing the integrand along its tangent and normal directions. For any $\boldsymbol{\varphi} \in \mathcal{D}(\Omega_{\mathcal{M}_+}^\epsilon, \mathbb{R}^3)$ one has (let us temporarily set $\tau_3(\sigma) := \nu(\sigma)$ to shorten notation)

$$\boldsymbol{h}_\epsilon[u] \cdot (\operatorname{\mathbf{curl}} \boldsymbol{\varphi} \circ \psi_\epsilon), \quad = \quad \frac{1}{2} \sum_{i \in \mathbb{N}_3} (\boldsymbol{h}_\epsilon[u] \times \tau_i) \cdot (\operatorname{\mathbf{curl}} \boldsymbol{\varphi} \circ \psi_\epsilon \times \tau_i). \tag{69}$$

We then denote by $\nabla_{\mathrm{skw}}$ the skew-symmetric part of the gradient defined by $\nabla_{\mathrm{skw}} \boldsymbol{\varphi} := (\nabla^\mathsf{T} \boldsymbol{\varphi} - \nabla \boldsymbol{\varphi})/2$. From (69) we obtain

$$\boldsymbol{h}_\epsilon[u] \cdot (\operatorname{\mathbf{curl}} \boldsymbol{\varphi} \circ \psi_\epsilon) \quad = \quad \sum_{i,j \in \mathbb{N}_3} (\boldsymbol{h}_\epsilon[u] \times \tau_i) \cdot ((\nabla_{\mathrm{skw}} \boldsymbol{\varphi} \circ \psi_\epsilon) \tau_i \cdot \tau_j) \tau_j \tag{70}$$

$$= \quad 2 \sum_{i<j \in \mathbb{N}_3} ((\boldsymbol{h}_\epsilon[u] \times \tau_i) \cdot \tau_j)((\nabla_{\mathrm{skw}} \boldsymbol{\varphi} \circ \psi_\epsilon) \tau_i \cdot \tau_j). \tag{71}$$

Next, we compute the relation between $\nabla_{\mathrm{skw}} \boldsymbol{\varphi} \circ \psi_\epsilon$ and $\nabla_\mathcal{M}(\boldsymbol{\varphi} \circ \psi_\epsilon)$. To this end, let us first note that for any $\boldsymbol{\varphi} \in \mathcal{D}(\mathcal{M}_+, \mathbb{R}^3)$, the function $\boldsymbol{\varphi}_\epsilon := \boldsymbol{\varphi} \circ \psi_\epsilon^{-1}$ is in $\mathcal{D}(\Omega_{\mathcal{M}_+}^\epsilon, \mathbb{R}^3)$, and moreover

$$(\nabla_{\mathrm{skw}} \boldsymbol{\varphi}_\epsilon \circ \psi_\epsilon) \tau_1 \cdot \tau_2 \quad = \quad (d\boldsymbol{\varphi}_\epsilon \circ \psi_\epsilon) \tau_1 \cdot \tau_2 - (d\boldsymbol{\varphi}_\epsilon \circ \psi_\epsilon) \tau_2 \cdot \tau_1$$

$$= \quad \frac{1}{1 + \epsilon s \kappa_1} \partial_{\tau_1}(\boldsymbol{\varphi} \circ \psi_\epsilon) \cdot \tau_2 - \frac{1}{1 + \epsilon s \kappa_2} \partial_{\tau_2}(\boldsymbol{\varphi} \circ \psi_\epsilon) \cdot \tau_1. \tag{72}$$

Similarly, we compute the normal component of the tangential image of $\nabla_{\mathrm{skw}}$. For any $i \in \mathbb{N}_2$ we get

$$(\nabla_{\mathrm{skw}} \boldsymbol{\varphi}_\epsilon \circ \psi_\epsilon) \tau_i \cdot \nu \quad = \quad (d\boldsymbol{\varphi}_\epsilon \circ \psi_\epsilon) \tau_i \cdot \nu - (d\boldsymbol{\varphi}_\epsilon \circ \psi_\epsilon) \nu \cdot \tau_i$$

$$= \quad \frac{1}{1 + \epsilon s \kappa_i} \partial_{\tau_2} \boldsymbol{\varphi} \cdot \nu - \frac{1}{\epsilon} \partial_s \boldsymbol{\varphi} \cdot \tau_i. \tag{73}$$

By taking the limit for $\epsilon \to 0$, of both members of (72) and (73), we get

$$\begin{aligned} \lim_{\epsilon \to 0} \left[ \epsilon \left( \nabla_{\mathrm{skw}} \boldsymbol{\varphi}_\epsilon \circ \psi_\epsilon \right) \tau_1 \cdot \tau_2 \right] &= 0, \\ \lim_{\epsilon \to 0} \left[ \epsilon (\nabla_{\mathrm{skw}} \boldsymbol{\varphi}_\epsilon \circ \psi_\epsilon) \tau_i \cdot \nu \right] &= -\partial_s \boldsymbol{\varphi} \cdot \tau_i, \end{aligned} \tag{74}$$

for any $i \in \mathbb{N}_2$. As $(\boldsymbol{h}_\epsilon[u_\epsilon])_{\epsilon \in I_\delta}$ satisfies (60), taking into account (71), we have that for every $\boldsymbol{\varphi} \in \mathcal{D}(\mathcal{M}_+, \mathbb{R}^3)$ (here we set as before $\boldsymbol{\varphi}_\epsilon := \boldsymbol{\varphi} \circ \psi_\epsilon$)

$$0 = \epsilon \int_{\mathcal{M}_+} \boldsymbol{h}_\epsilon[u_\epsilon] \cdot (\operatorname{\mathbf{curl}} \boldsymbol{\varphi}_\epsilon \circ \psi_\epsilon) \, \mathfrak{g}_\epsilon \, d\sigma \, ds \tag{75}$$

$$= 2 \sum_{i<j \in \mathbb{N}_3} \int_{\mathcal{M}_+} [(\boldsymbol{h}_\epsilon[u_\epsilon] \times \tau_i) \cdot \tau_j] \left[ \epsilon (\nabla_{\mathrm{skw}} \boldsymbol{\varphi}_\epsilon \circ \psi_\epsilon) \tau_i \cdot \tau_j \right] \mathfrak{g}_\epsilon \, d\sigma \, ds. \tag{76}$$

Since $\boldsymbol{h}_\epsilon[u_\epsilon] \rightharpoonup \boldsymbol{h}_0$ weakly in $L^2(\mathbb{R}^3, \mathbb{R}^3)$, taking into account (41) and (74) and passing to the limit for $\epsilon \to 0$ in the previous expression, we finish with the relation

$$\int_{\mathcal{M}_+} (\boldsymbol{h}_0(\sigma,s) \times \nu(\sigma)) \cdot \partial_s \boldsymbol{\varphi}(\sigma,s) \, \mathrm{d}\sigma \, \mathrm{d}s = 0 \quad \forall \boldsymbol{\varphi} \in \mathcal{D}(\mathcal{M}_+, \mathbb{R}^3), \tag{77}$$

from which we deduce that the quantity $\boldsymbol{h}_0 \times \nu$ does not depend on the $s$-variable. Since $\boldsymbol{h}_0 \times \nu \in L^2(\mathcal{M}_+)$ one necessarily has $\boldsymbol{h}_0 \times \nu = 0$. Hence, the weak limit $\boldsymbol{h}_0$ has no tangential component and that means (cf. (18)) that $\lim_{\epsilon \to 0} \mathcal{W}_\tau^\epsilon(u_\epsilon) = 0$. We have so proved the following result.

LEMMA 10. *If $u_\epsilon \rightharpoonup u_0$ weakly in $H^1(\mathcal{M}, \mathbb{S}^2)$ and $\liminf_{\epsilon \to 0} \mathcal{V}_\epsilon(u_\epsilon) < +\infty$, then $\lim_{\epsilon \to 0} \mathcal{W}_\tau^\epsilon(u_\epsilon) = 0$ and*

$$\mathcal{W}_0 = \lim_{\epsilon \to 0} \mathcal{W}_\mathcal{M}^\epsilon(u_\epsilon) = \lim_{\epsilon \to 0} \mathcal{W}_\nu^\epsilon(u_\epsilon) = \int_S (u_0(\sigma) \cdot \nu(\sigma))^2 \, \mathrm{d}\sigma. \tag{78}$$

### 4.3. The Γ-limit of the family $(\mathcal{F}_\epsilon)_{\epsilon \in I_\delta}$

We complete the proof of Theorem 1 by showing that $\mathcal{F}_0$ is given by (24). As pointed out in Remark 7, it is sufficient to show that Γ-$\lim \mathcal{V}_\epsilon = \mathcal{E}_0 + \mathcal{W}_0$. We note that if $u_\epsilon \rightharpoonup u$ in $H^1(\mathcal{M}, \mathbb{S}^2)$ and $\partial_s u \neq 0$ then $\liminf_{\epsilon \to 0} \mathcal{V}_\epsilon(u_\epsilon) = +\infty$ and therefore the Γ-liminf inequality is trivially satisfied. On the other hand, if $\partial_s u = 0$, then from Lemma 8 and Lemma 10 we get

$$\liminf_{\epsilon \to 0} \mathcal{V}_\epsilon(u_\epsilon) = \liminf_{\epsilon \to 0} (\mathcal{E}_\tau^\epsilon(u_\epsilon) + \mathcal{E}_\nu^\epsilon(u_\epsilon)) + \lim_{\epsilon \to 0} (\mathcal{W}_\nu^\epsilon(u_\epsilon) + \mathcal{W}_\tau^\epsilon(u_\epsilon)) \tag{79}$$

$$\geqslant \mathcal{E}_0(u_0) + \mathcal{W}_0(u_0) \tag{80}$$

$$= \mathcal{V}_0(u_0).$$

Finally, for any $u_0 \in H^1(\mathcal{M}, \mathbb{S}^2)$ such that $\partial_s u_0 = 0$, the constant (with respect to the index $\epsilon$) family $(u_\epsilon)_{\epsilon \in I_\delta} = (u_0)_{\epsilon \in I_\delta}$ is a recovery sequence. Indeed we have

$$\limsup_{\epsilon \to 0} \mathcal{V}_\epsilon(u_0) = \limsup_{\epsilon \to 0} (\mathcal{E}_\tau^\epsilon(u_0) + \mathcal{E}_\nu^\epsilon(u_0)) + \lim_{\epsilon \to 0} (\mathcal{W}_\nu^\epsilon(u_\epsilon) + \mathcal{W}_\tau^\epsilon(u_\epsilon)) = \mathcal{V}_0(u_0), \tag{81}$$

and this completes the proof of Theorem 1.

### 5. Conclusion and Acknowledgment

We have computed the Γ-limit of the micromagnetic energy functional when the shell is generated by a bounded and smooth convex surface. Our result provides a solid ground to most of the studies on nanomagnets with curved shape which are currently under investigation by the theoretical physics community (e.g., [27, 29, 37, 36, 39, 42, 43, 48]; this list is certainly far from complete).

In particular, our result validates the variational model widely used in the analysis of magnetic thin spherical shells which are currently worth of interest due to their capability to support skyrmion solutions (see e.g. [37]). Indeed, for a magnetic spherical film with perpendicular magnetocrystalline anisotropy, the Γ-limit (after a suitable rescaling and in the absence of an external applied field) reads as

$$\mathcal{F}_0 : u \in H^1(\mathbb{S}^2, \mathbb{S}^2) \mapsto \int_{\mathbb{S}^2} |\nabla_\tau u|^2 + \kappa^2 (u \cdot \nu)^2 \, \mathrm{d}\sigma \tag{82}$$

with $\kappa^2$ summarizing the contributions of both the crystal and shape anisotropy. The investigation of the metastable states of (82) turns out to be a challenging problem with far-reaching consequences in the modern magnetic storage technology. Indeed, as some formal asymptotics shows, $\mathcal{F}_0$ exhibits two topologically protected metastable states, known as the *vortex* and the *onion* state, depending on the value of $\kappa^2$ ($\kappa^2_{\text{onion}} < \kappa^2_{\text{vortex}}$, cf. Figure 3). These states are characterized by distinct skyrmion numbers [46] and therefore appropriate for the design of future racetrack memory devices [26].

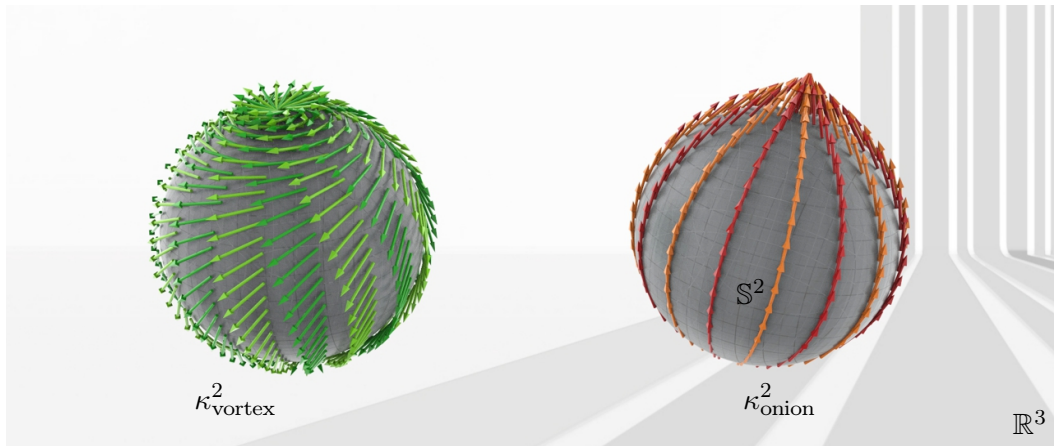

**Figure 3.** Magnetic spherical films are currently worth of interest due to their capability to support skyrmion solutions. (From left to right) The *vortex* and the *onion* state ($\kappa^2_{\text{onion}} < \kappa^2_{\text{vortex}}$).

The formation of these two states can be heuristically explained as follows. Let us recall that a ferromagnetic particle occupying a spherical region can support constant in space magnetizations. In other words, if the ferromagnet occupies the spherical region $B_r$ (of radius $r$) and $m$ is constant in $B_r$, then the induced demagnetizing field $\boldsymbol{h}_{\text{d}}[m]$ is also constant in $B_r$ (see [21, 33, 41]). Moreover, according to Brown's fundamental theorem of the theory of fine ferromagnetic particles (cf. [3, 4, 6, 13, 22]), there exists a critical radius $r_c$ below which both the global and the local minimizers of the micromagnetic energy functional $\mathcal{G}(\cdot, B_r)$ are constant in space. Also, there exists a critical radius $R_c > r_c$ such that when $r > R_c$ the global minimizers of $\mathcal{G}(\cdot, B_r)$ are no more constant in space and *vortex* type solutions start to be energetically preferable (cf. [13]).

Now, for a spherical shell $\Omega_{\alpha,r} := B_r \setminus B_{\alpha r}$, with $\alpha < 1$, one has

$$\mathcal{G}(\cdot, \Omega_{\lambda,r}) = \mathcal{G}(\cdot, B_r) - \mathcal{G}(\cdot, B_{\alpha r}), \tag{83}$$

and the relations $r < r_c$ and $r > R_c$ translate, in the framework of the limiting energy functional (82), as $\kappa^2 < \kappa_c^2$ and $\kappa^2 > K_c^2$ for some suitable constants $\kappa_c^2$ and $K_c^2$. Therefore, when $r < r_c$, i.e., when $\kappa^2$ is sufficiently small, the ground states are the result of an energetic competition among constant in space magnetizations, which tend to minimize $\mathcal{G}(\cdot, B_r)$, and normal (to the sphere) configurations which tend to maximize $\mathcal{G}(\cdot, B_{\alpha r})$. This leads to *onion* type configurations. On the other hand, when $r > R_c$, i.e., when $\kappa^2$ is sufficiently large, the minimizers are the result of an energetic competition among vortex type configurations, which tend to minimize $\mathcal{G}(\cdot, B_r)$, and the normal ones which tend to maximize $\mathcal{G}(\cdot, B_{\alpha r})$. This leads to the *vortex* type configurations. This and many other aspects of the question will be the object of forthcoming works.

## 5.1. Acknowledgment

This work was partially supported by the labex LMH through the grant no. ANR-11-LABX-0056LMH in the *Programme des Investissements d'Avenir*. The author gratefully acknowledges the support provided by Severo Ochoa Program SEV-2013-0323 and Basque Government BERC Program 2014-2017.

## Bibliography

[1] E. ACERBI, I. FONSECA, and G. MINGIONE, *Existence and regularity for mixtures of micromagnetic materials*, Proc. R. Soc., 462 (2006), pp. 2225–2243.

[2] M. S. AGRANOVICH, *Sobolev Spaces, Their Generalizations and Elliptic Problems in Smooth and Lipschitz Domains*, Springer, Cham, 2015.

[3] A. AHARONI, *Elongated single-domain ferromagnetic particles*, Journal of applied physics, 63 (1988), pp. 5879–5882.

[4] F. ALOUGES and K. BEAUCHARD, *Magnetization switching on small ferromagnetic ellipsoidal samples*, ESAIM: Control, Optimisation and Calculus of Variations, 15 (2009), pp. 676–711.

[5] F. ALOUGES and G. DI FRATTA, *Homogenization of composite ferromagnetic materials*, Proc. R. Soc. A, 471 (2015).

[6] F. ALOUGES, G. DI FRATTA, and B. MERLET, *Liouville type results for local minimizers of the micromagnetic energy*, Calculus of Variations and Partial Differential Equations, 53 (2015), pp. 525–560.

[7] F. ALOUGES and S. LABBÉ, *Convergence of a ferromagnetic film model*, Comptes Rendus Mathematique, 344 (2007), pp. 77–82.

[8] G. BERTOTTI, *Hysteresis in magnetism: for physicists, materials scientists, and engineers*, Academic Press, San Diego, 1998.

[9] G. BERTOTTI, I. MAYERGOYZ, and C. SERPICO, *Nonlinear magnetization dynamics in nanosystems*, Elsevier Science Limited, Amsterdam, 2009.

[10] A. BRAIDES and A. DEFRANCESCHI, *Homogenization of multiple integrals*, vol. 263, Clarendon Press, Oxford, 1998.

[11] W. F. BROWN, *Magnetostatic Principles in Ferromagnetism*, North-Holland, Amsterdam, 1962.

[12] W. F. BROWN, *Micromagnetics*, Interscience Publishers, London, 1963.

[13] W. F. BROWN, *The fundamental theorem of the theory of fine ferromagnetic particles*, Annals of the New York Academy of Sciences, 147 (1969), pp. 463–488.

[14] G. CARBOU, *Thin layers in micromagnetism*, Mathematical Models and Methods in Applied Sciences, 11 (2001), pp. 1529–1546.

[15] G. DAL MASO, *Introduction to $\Gamma$-convergence*, vol. 8, Birkhauser, Boston, 1993.

[16] R. DAUTRAY and J. LIONS, *Mathematical Analysis and Numerical Methods for Science and Technology: Volume 3 Spectral Theory and Applications*, Mathematical Analysis and Numerical Methods for Science and Technology, Springer Berlin Heidelberg, 1999.

[17] R. DAUTRAY and J. LIONS, *Mathematical Analysis and Numerical Methods for Science and Technology: Volume 4 Integral Equations and Numerical Methods*, Mathematical Analysis and Numerical Methods for Science and Technology, Springer Berlin Heidelberg, 1999.

[18] A. DESIMONE, *Hysteresis and imperfection sensitivity in small ferromagnetic particles*, Meccanica, 30 (1995), pp. 591–603.

[19] A. DESIMONE, R. V. KOHN, S. MÜLLER, and F. OTTO, *A reduced theory for thin-film micromagnetics*, Communications on pure and applied mathematics, 55 (2002), pp. 1408–1460.

[20] A. DESIMONE, R. V. KOHN, S. MÜLLER, and F. OTTO, *Recent analytical developments in micromagnetics*, The science of hysteresis, 2 (2006), pp. 269–381.

[21] G. DI FRATTA, *The newtonian potential and the demagnetizing factors of the general ellipsoid*, Proceedings of the Royal Society of London A: Mathematical, Physical and Engineering Sciences, 472 (2016).

[22] G. DI FRATTA, C. SERPICO, and M. D'AQUINO, *A generalization of the fundamental theorem of Brown for fine ferromagnetic particles*, Physica B: Condensed Matter, 407 (2012), pp. 1368–1371.

[23] M. DO CARMO, *Differential geometry of curves and surfaces*, Prentice-hall, Englewood Cliffs, 1976.

[24] I. DZYALOSHINSKY, *A thermodynamic theory of "weak" ferromagnetism of antiferromagnetics*, Journal of Physics and Chemistry of Solids, 4 (1958), pp. 241–255.

[25] H. FEDERER, *Geometric Measure Theory*, Classics in Mathematics, Springer Berlin Heidelberg, 2014.
[26] A. FERT, V. CROS, and J. SAMPAIO, *Skyrmions on the track*, Nature nanotechnology, 8 (2013), pp. 152–156.
[27] Y. GAIDIDEI, V. P. KRAVCHUK, and D. D. SHEKA, *Curvature effects in thin magnetic shells*, Physical review letters, 112 (2014), p. 257203.
[28] G. GIOIA and R. D. JAMES, *Micromagnetics of very thin films*, Proceedings of the Royal Society of London A: Mathematical, Physical and Engineering Sciences, 453 (1997), pp. 213–223.
[29] A. GOUSSEV, J. M. ROBBINS, and V. SLASTIKOV, *Domain wall motion in thin ferromagnetic nanotubes: analytic results*, EPL (Europhysics Letters), 105 (2014), p. 67006.
[30] A. HUBERT and R. SCHÄFER, *Magnetic domains: the analysis of magnetic microstructures*, Springer Verlag, New York, 2008.
[31] R. IGNAT, *A survey of some new results in ferromagnetic thin films*, in Séminaire: Équations aux dérivées partielles, 2007-2008 (2009), pp. 1–19.
[32] J. JACKSON, *Classical Electrodynamics*, Wiley, 2007.
[33] O. KELLOGG, *Foundations of potential theory*, vol. 31, Dover Publications, New York, 2010.
[34] W. KLINGENBERG, *A course in differential geometry*, vol. 51, Springer-Verlag, New York, 2013.
[35] R. V. KOHN and V. SLASTIKOV, *Another thin-film limit of micromagnetics*, Archive for rational mechanics and analysis, 178 (2005), pp. 227–245.
[36] V. P. KRAVCHUK, D. D. SHEKA, R. STREUBEL, D. MAKAROV, O. G. SCHMIDT, and Y. GAIDIDEI, *Out-of-surface vortices in spherical shells*, Physical Review B, 85 (2012), p. 144433.
[37] V. P. KRAVCHUK, U. K. RÖSSLER, O. M. VOLKOV, D. D. SHEKA, J. VAN DEN BRINK, D. MAKAROV, H. FUCHS, H. FANGOHR, and Y. GAIDIDEI, *Topologically stable magnetization states on a spherical shell: curvature-stabilized skyrmions*, Physical Review B, 94 (2016), p. 144402.
[38] L. LANDAU and E. LIFSHITZ, *On the theory of the dispersion of magnetic permeability in ferromagnetic bodies*, Phys. Z. Sowjetunion, 8 (1935), pp. 101–114.
[39] P. LANDEROS, S. ALLENDE, J. ESCRIG, E. SALCEDO, D. ALTBIR, and E. VOGEL, *Reversal modes in magnetic nanotubes*, Applied Physics Letters, 90 (2007), p. 102501.
[40] T. MORIYA, *Anisotropic superexchange interaction and weak ferromagnetism*, Physical Review, 120 (1960), p. 91.
[41] J. A. OSBORN, *Demagnetizing factors of the general ellipsoid*, Physical Review, 67 (1945), pp. 351–357.
[42] D. D. SHEKA, V. P. KRAVCHUK, and Y. GAIDIDEI, *Curvature effects in statics and dynamics of low dimensional magnets*, Journal of Physics A: Mathematical and Theoretical, 48 (2015), p. 125202.
[43] D. D. SHEKA, V. P. KRAVCHUK, M. I. SLOIKA, and Y. GAIDIDEI, *Equilibrium states of soft magnetic hemispherical shell*, Spin, 3 (2013), p. 1340003.
[44] V. SLASTIKOV, *Micromagnetics of thin shells*, Mathematical Models and Methods in Applied Sciences, 15 (2005), pp. 1469–1487.
[45] V. SLASTIKOV and C. SONNENBERG, *Reduced models for ferromagnetic nanowires*, IMA Journal of Applied Mathematics, 77 (2012), pp. 220–235.
[46] M. SLOIKA, D. SHEKA, V. KRAVCHUK, O. PYLYPOVSKYI, and Y. GAIDIDEI, *Geometry induced phase transitions in magnetic spherical shell*, Journal of Magnetism and Magnetic Materials, 443 (2017), pp. 404–412.
[47] A. VISINTIN, *On Landau-Lifshitz'equations for ferromagnetism*, Japan Journal of Industrial and Applied Mathematics, 2 (1985), pp. 69–84.
[48] M. YAN, C. ANDREAS, A. KÁKAY, F. GARCA-SÁNCHEZ, and R. HERTEL, *Fast domain wall dynamics in magnetic nanotubes: suppression of walker breakdown and cherenkov-like spin wave emission*, Applied physics letters, 99 (2011), p. 122505.